\documentclass[dvips, 12pt,a4paper]{article}

\usepackage[greek,american]{babel}
\usepackage{amssymb}
\usepackage{amsfonts}
\usepackage{amsmath}
\usepackage{subfigure}
\usepackage[dvips]{graphicx}

\def\mP{\mathbb P}

\def\E{\mathbb E}

\begin{document}

\begin{center}
{{{\Large\sf\bf Dependent Random Density Functions with Common Atoms and Pairwise Dependence}}}\\
\vspace{0.5cm}
{\large\sf Spyridon J. Hatjispyros $^*$, Theodoros Nicoleris$^{**}$ and Stephen G. Walker$^{***}$}\\
\end{center}

\centerline{\sf $^{*}$ Department of Mathematics, University of the Aegean,}
\centerline {\sf Karlovassi, Samos, GR-832 00, Greece.} 
\centerline{\sf $^{**}$ Department of Economics, National and Kapodistrian University of Athens,}
\centerline{\sf Athens, GR-105 59, Greece. }  
\centerline{\sf $^{***}$Department of Mathematics, University of Texas at Austin,}
\centerline{\sf Austin, Texas 7812, USA. }

\begin{abstract} 
The paper is concerned with constructing pairwise dependence between $m$ random density functions each of which is modeled as a mixture of Dirichlet process model.  The key to this is how to create dependencies between random Dirichlet processes. The present paper adopts a plan previously used for creating pairwise dependence, with the simplification that all random Dirichlet processes share the same atoms. Our contention is that for all dependent Dirichlet process models, common atoms are sufficient. 

We show that by adopting common atoms, it is possible to compute the $L_p$ distances between all pairs of random probability measures.   

\vspace{0.1in} \noindent {\sl Keywords:} Bayesian nonparametric inference; Dependent Dirichlet process; $L_p$ distance; Mixture of Dirichlet process; Pairwise dependence.
\end{abstract}

\vspace{0.5in} \noindent {\bf 1. Introduction.} There has been substantial recent interest in the construction of dependent probability measures.  We first restrict our thoughts to dependence between two random probability measures, $\mP_1$ and $\mP_2$, though for the paper we are concerned with creating (pairwise) dependence between $m$ probability measures. To complete the motivation, the typical scenario in which such measures are employed is with mixture models, generating random densities $f_1$ and $f_2$, whereby
$$
f_1(x)=\int_\Theta K(x|\theta)\,\mP_1(d \theta)\quad{\rm and}\quad 
f_2(x)=\int_\Theta K(x|\theta)\,\mP_2(d \theta)
$$
The marginal models for  each $f_j$ are a random density function based on the benchmark mixture of Dirichlet process
model (Lo, 1984); so that each
$\mP_j$ is a Dirichlet process (Ferguson, 1973), and $K(x|\theta)$ is a density function for each $\theta\in\Theta$.
The reason for dependence is that it is thought that properties of $f_1$ and $f_2$ are similar in some way; for example the means are similar, or it is thought the distance between them is small, in that they resemble each other.

We can write each ${\mathbb P}_{j}$ using the constructive definition of the Dirichlet
process given in Sethuraman (1994); so that  
$$
{\mathbb P}_{j}=\sum_{k=1}^\infty w_{jk}\,\delta_{\theta_{jk}}
$$
where we write $\theta_{j}=(\theta_{jk})_{k=1}^\infty$, being independent
and identically distributed from some fixed distribution
$P_0(\theta)$, with density function $p_0(\theta)$. 
And we write ${\bf w_{j}}=(w_{jk})_{k=1}^\infty$; a stick--breaking process,
so if all the $(z_{jk})_{k=1}^\infty$ are independent and identically distributed from the beta$(1,c)$ distribution, for some $c>0$, then $w_{j1}=z_{j1}$
and, for $k>1$,
$$
w_{jk}=z_{jk}\prod_{j<l}(1-z_{jl}).
$$
In applications, the $\mP_1$ and $\mP_2$ can be made dependent in a variety of ways. The modeling of dependent random Dirichlet processes has been the focus of much recent research in Bayesian nonparametrics; following original developments in MacEachern (1999).  More recent work is to be found in De Iorio et al. (2004), Griffin and Steel (2006), Dunson and Park (2008) and summaries in Hjort et al. (2010). The use of the Dependent Dirichlet process arises mostly in regression problems
where a random measure $\mP_z$ is constructed for each covariate $z$ and is of the type
$$
{\mathbb P}_{z}=\sum_{k=1}^\infty w_{k}(z)\,\delta_{\theta_{k}(z)},
$$   
where $(w_k(z),\theta_k(z))$ are processes of weights and atoms.

On the other hand, the work reported here is about modeling a finite number of densities, equivalent to the regression model with a finite number of fixed covariates.
Previous work related to this type of structure is to be found in M\"uller et al. (2004) and Bulla et al. (2009), and more recently Hatjispyros et al. (2011), Kolossiatis et al. (2013) and Griffin et al. (2013).
These papers model an arbitrary but finite number of random distribution functions, $(\mP_1,\ldots,\mP_m)$
via a common component and an index specific idiosyncratic component. That is, the basic idea is to model
$$\mP_j=p_j\mP_0+(1-p_j)\mP_j^*$$  
 where $\mP_0$ is the common component to all other distributions and $(\mP_j^*)$ are the idiosyncratic parts, unique to each $\mP_j$, $p_j \in (0,1) $. And $(\mP_0,\mP_j^*)$ are modeled as mutually independent mixture of Dirichlet process models, or based on some other stick--breaking process (Ishwaran and James, 2001) or normalized random measure as in the case of Griffin et al. (2013).
In all these cases the (random) atoms for the underlying Dirichlet process are unique to each $j$. 

The present paper is concerned with constructing $\mP_j$ so that there is a unique common component for each pair
$(\mP_j,\mP_{j'})$ with $j\ne j'$. This allows arbitrary pairwise dependence between any two $\mP_j$ and $\mP_{j'}$. 
The details of this idea for the case $m=2$  were presented in Hatjispyros et al. (2011). Our  numeric results target the key $m=3$ case, which essentially gives
access to arbitrary values of $m$. It is clear that $m=2$ is a straightforward special case.

Moreover, our thesis is that when constructing e.g. $\mP_1$ and $\mP_2$, it is sufficient to create a dependence between the weights and to take $\theta_1=\theta_2$. That is, it is sufficient for the random probability measures $\mP_1$ and $\mP_2$ to have the same atoms. We demonstrate this sufficiency throughout the paper.

The key to the understanding of this idea is quite straightforward. From a fixed set of random atoms obtained for one probability measure, another probability measure can be obtained by reassigning the weights to these same atoms. 
It should be clear that varying weights can provide probability measures which are either remarkably close, when weights are similar, to probability measures which are far apart, when the weights are dissimilar. 

It is less interpretable to have both varying weights and atoms. For if the weights are similar, there is nothing to be said about the closeness of the distributions as the ulimate picture will depend on the atoms. And we decide to make the atoms similar and allow the weights to vary so that with the same atoms we can easily compute distances between the two measures. We will look at this in Section 2.

We will demonstrate the idea using a dependent model suggested by  Hatjispyros et al. (2011) for constructing pairwise dependence for the finite set of densities
$$(f_1,\ldots,f_m)$$
where each $f_j$ is a random density function. The idea is that from each density we observe independent data sets, yet the densities from which they come share common features. For example, they may all have similar tail behaviour or even have common variances, and so on. Hence, it is imperative to model in the prior an arbitrary level of dependence between each pair $(f_j,f_k)$, for each $j\ne k$. The key is to construct the prior model in such a way that for every pair there is a ``common" part and a ``difference" part, to be explained more explicitly in due course.

The layout of the paper is as follows. In Section 2 we provide some preliminary findings concerned with the evaluation of distances between probability measures being generated. This is possible when atoms are common to each distribution.
In section 3 we describe the model of interest and introduce the latent variables which will form the basis of the Gibbs sampler which is described in Section 4. Section 5 contains a numerical illustration involving a real data set and finally Section 6 concludes with a summary and future work.

\vspace{0.2in} \noindent
{\bf 2. Preliminaries.} Given the common atoms for $\mP_1$ and $\mP_2$ we can easily compute distances between them and also between the corresponding mixed density functions. So, for suitable sets $A$,
$$|\mP_1(A)-\mP_2(A)|=\left|\sum_{\theta_j\in A}(w_{1j}-w_{2j})\right|.$$
Therefore, for the total variation distance between $\mP_1$ and $\mP_2$ we have
$$
d(\mP_1,\mP_2)=\sup_A|\mP_1(A)-\mP_2(A)|=\sup_J|W_{1J}-W_{2J}|,
$$
where, for example, 
$$
W_{1J}=\sum_{j\in J}w_{1j},
$$
where $J$ is an index set.
This is simple to interpret but impossible to obtain and control if the atoms are not identical. In fact, if the atoms are disjoint for each measure, close weights, even identical weights, says nothing about how close the two measures are to each other. We believe equal atoms is fundamental as a consequence. 

\vspace{0.2in}
We can formalize the idea of using common atoms and the sufficiency of this by considering 
the following

\smallskip\noindent
{\bf Lemma $1$:~}{\sl We consider the densities $f_0(x)=\sum_{j=1}^\infty q_j\,K(x|\theta_j)$ and
$f_1(x)=\sum_{j=1}^\infty w_j\,K(x|\psi_j)$
where the $(\theta_j)$, $(\psi_j)$ and $(w_j)$ are fixed, and the $(\theta_j)$ are dense in $\Theta$.
Then we can find $(q_j)$ such that the $L_1$-distance between $f_0$ and $f_1$ can be made arbitrarily small.}

\smallskip\noindent
{\bf Proof:~} Because the $(\theta_j)$ are dense in $\Theta$,  
for any $j$ and any $\epsilon_j>0$, we can find $l(j)$ and $\theta_{l(j)}$  such that
$$
d_1\big(K(\,\cdot\,|\psi_j),K(\,\cdot\,|\theta_{l(j)})\big)<\epsilon_j.
$$

We assume without loss of generality  that $l(j)$ is one to one, for $j\in {\mathbb N}$. 
Then, we put $q_{l(j)}=w_j$,  and for the remaining  indices  we set  $q_k=0$. 
Thus, we have that $ \sum_{j=1}^\infty q_{l(j)}=\sum_{j=1}^\infty w_j=1$. Hence, 
$$
\begin{array}{ll}
d_1(f_0,f_1)  & ={\large \int}_X \left|\sum_{j=1}^\infty w_j \big\{K(x|\psi_j)-K(x|\theta_{l(j)})\big\}\right|\,dx \\ \\
& \leq \sum_{j=1}^\infty w_j\,d_1\big(K(\,\cdot\,|\psi_j),K(\,\cdot\,|\theta_{l(j)})\big)
  \leq \sum_{j=1}^\infty w_j\,\epsilon_j.
\end{array}
$$
Now choose $(\epsilon_j)$ such that, for any $\epsilon>0$, we have 
$$
\sum_{j=1}^\infty w_j\epsilon_j<\epsilon,
$$
and the lemma follows.

\smallskip\noindent
Thus, even though atoms are fixed across densities, as long as they form a dense set, we can approximate arbitrarily accurately any density with any atoms.

Hence, we can obtain weights to allow this distance to be either 0, the weights coincide, or to be 1. A dependent prior for $(\bf{w_1,w_2})$  can be used to provide a small distance between $\mP_1$ and $\mP_2$ if that is what is required. This is not so obvious to achieve if the atoms are dissimilar. Moreover, if atoms are dissimilar computing the informative distances is not possible and one is left with computing objects such as
$$
\mbox{Cov}(\mP_1(A),\mP_2(A)),
$$
which, although  it provides some insight on the dependence between  $\mP_1$  and $\mP_2$,  is not the appropriate learning tool for the   similarity between two random densities. Clearly, in that case the total variation distance
(or the $L_2$ distance between densities)  is a  more  appropriate measure.   
 
We now turn to looking at a distance between $f_1$ and $f_2$, and for ease of computation we make this the $L_2$ distance,
$$
d_2(f_1,f_2)=\int_X \left(f_1(x)-f_2(x) \right)^2 d x.
$$
\smallskip\noindent
{\bf Lemma $2$:~}{\sl We consider the random mixtures $f_i(x)=\sum_{j=1}^\infty w_{ij}K(x|\theta_j),\,i=1,2$ with 
$\theta_j$ independent, coming from the base measure $P_0$ for all $k\ge 1$, then
$$
{\E} \left[\, d_2(f_1,f_2)\, | \bf{w_1,w_2}\right]\,=\,
(\alpha-\beta)\sum_{j=1}^\infty (w_{1j}-w_{2j})^2, 
$$
where $\alpha=\E\left\{\int_X K(x|\theta_j)^2 dx\right\}$ and 
$\beta=\E\left\{\int_X K(x|\theta_j)K(x|\theta_k)dx\right\}$.}

\smallskip\noindent
{\bf Proof:~} It is that
$$
\begin{array}{ll}
d_2(f_1,f_2) &  =\int_X \left\{\sum_{j=1}^\infty (w_{1j}-w_{2j})K(x|\theta_j) \right\}^2 dx \\ \\
& = \sum_{j=1}^\infty \Phi_j^2 \,(w_{1j}-w_{2j})^2+
\sum_{j\ne k}\Phi_{jk}\,(w_{1j}-w_{2j})(w_{1k}-w_{2k}),
\end{array}
$$ 
where
$$
\Phi_j^2=\int_X K(x|\theta_j)^2 dx\quad\mbox{and}\quad 
\Phi_{jk}=\int_X K(x|\theta_j)K(x|\theta_k) dx.
$$
Taking expectations over the atoms but keeping the weights fixed, we have
$$
\begin{array}{ll}
\E\left[\, d_2(f_1,f_2)\, |{\bf w_1, w_2} \right]\, & =\,
\alpha\sum_{j=1}^\infty(w_{1j}-w_{2j})^2+\beta\sum_{j\ne k}(w_{1j}-w_{2j})(w_{1k}-w_{2k})\\ \\
& =\,
(\alpha-\beta)\sum_{j=1}^\infty (w_{1j}-w_{2j})^2+\beta{\cal S},
\end{array}
$$
where 
$$
{\cal S}=\sum_{j=1}^\infty(w_{1j}-w_{2j})^2+2\sum_{j<k}(w_{1j}-w_{2j})(w_{1k}-w_{2k}),\quad
\alpha=\E\left(\Phi_j^2\right)\quad{\rm and}\quad 
\beta={\mathbb E}\left(\Phi_{jk}\right),
$$ 
and note that $\alpha-\beta=\int_X\rm Var\left\{K(x|\theta_j)\right\}dx>0$.
Now it is very easy to show that ${\cal S}=0$ by
utilizing the identities 
$$
\sum_{j=1}^\infty w_{ij}^2+2\sum_{j<k}w_{ij}w_{ik}=1,
$$
and
$$
\sum_{j=1}^\infty w_{1j}w_{2j}+\sum_{j<k}w_{1j}w_{2k}+\sum_{j<k}w_{2j}w_{1k}=1,
$$
and the lemma follows.
 
So, again, crucial distances can be understood solely through the weights.

If we are now interested in creating dependent weights for $f_1$ and $f_2$ (which will allow extensions to a larger number of densities with pairwise dependence), then we construct weights of the type
$$
w_{1j}=pw_{11j}+(1-p)w_{12j}\quad\mbox{and}\quad w_{2j}=qw_{21j}+(1-q)w_{22j},
$$
with $w_{12j}=w_{21j}$. This gives a common part to $f_1$ and $f_2$, via a common part to $\mP_1$ and $\mP_2$, and we can write
$$
\mP_1=\sum_{j}w_{1j}\delta_{\theta_j}\quad\mbox{and}\quad\mP_2=\sum_{j}w_{2j}\delta_{\theta_j},
$$
due to the common atoms.

On the other hand, in  Hatjispyros et al. (2011) the model was described as
$$
\mP_1\,=\,p\,\mP_{11}+(1-p)\,\mP_{12}\quad\mbox{and}\quad 
\mP_2\,=\,q\,\mP_{21}+(1-q)\,\mP_{22}
$$
with $\mP_{12}=\mP_{21}$. It had to be like this due to the uncommon atoms of $\mP_{11}$, $\mP_{12}$ and $\mP_{22}$,
i.e.
$$\mP_{jl}=\sum_{k=1}^\infty w_{jlk}\delta_{\theta_{jlk}}.$$
The details of the model where random probability measures share common atoms is now described in Section 3.

\vspace{0.2in} \noindent
{\bf 3. The model.} We start off by describing the model as it was in Hatjispyros et al. (2011) and then proceed to detail the simplifications when atoms are common. So, we have the set of random density functions generated via
$$
f_j(x|\,{\mathbb P}_{j1},\ldots,{\mathbb P}_{j\,m} )\, = \,
\sum_{l=1}^m p_{jl}\,g_{j\,l}(x|\,{\mathbb P}_{j\,l}),\quad 1\le j\le m,
$$
where $\sum_{l=1}^m p_{jl}=1$. The random densities satisfy
$g_{jl}=g_{lj}$ and are independent mixtures of Dirichlet process models; so that
$$
g_{jl}(x)\,=\,g(x|\,{\mathbb P}_{jl})=\int_\Theta K(x|\,\theta)\,{\mathbb P}_{jl}(d\theta),
$$
for some kernel density $K(\cdot|\cdot)$ and $\{{\mathbb P}_{jl}:\,1\le j,\,l\le m\}$ form a matrix of random distributions with 
${\mathbb P}_{jl}={\mathbb P}_{lj}$ for $j>l$ and each other element is an independent Dirichlet process.
Equivalently, the random densities $(f_j)$ are dependent mixtures of the
dependent random measures
\begin{equation}
\label{Qmeasures}
{\mathbb Q}_j(d\theta) \, = \, \sum_{l=1}^m p_{jl}\,{\mathbb P}_{jl}(d\theta),\quad 1\le j\le m.
\end{equation}
In matrix notation
$$
{\mathbb Q}\,=\,{\bf A}\,{\bf 1},\quad {\bf A}\,=\,{\cal W}\otimes {\cal P}\,,
$$
where ${\cal W}=(p_{jl})$ is the matrix of weights and ${\cal P}=({\mathbb P}_{jl})$
is the symmetric matrix of the independent Dirichlet measures.

We can write each ${\mathbb P}_{jl}$ using the constructive definition of the Dirichlet
process given in Sethuraman (1994); so that, and now adopting common atoms for all probability measures,
$$
{\mathbb P}_{jl}=\sum_{k=1}^\infty w_{jlk}\,\delta_{\theta_{k}}
$$
where we write $\vartheta=(\theta_{k})_{k=1}^\infty$, being independent
and identically distributed from some fixed distribution
$P_0(\theta)$, with density function $p_0(\theta)$, and we write ${\bf w_{jl}}=(w_{jlk})_{k=1}^\infty$;
being a stick-breaking process;
so if all the $(z_{jlk})$ are independent and identically distributed from the beta$(1,c)$ distribution, for some $c>0$, then $w_{jl1}=z_{jl1}$
and, for $k>1$,
$$
w_{jlk}=z_{jlk}\prod_{r<k}(1-z_{jlr}).
$$
Hence, we can write  
$$
g_{jl}(x|\,{\bf w_{jl},\vartheta})\,=\,\sum_{k=1}^\infty w_{jlk}\,K(x|\,\theta_{k}),                             
$$
and
$$
f_j(x|\,({\bf w_{j1},\ldots, w_{jm})},\vartheta)\,=\,\sum_{k=1}^\infty
\left\{\sum_{l=1}^m p_{jl}w_{jlk}\right\}K(x|\,\theta_{k}).
$$
We could write this as
$$f_j(x|\,({\bf w_{j1},\ldots, w_{jm})},\vartheta)\,=\,\sum_{k=1}^\infty w_{jk}\,K(x|\,\theta_{k}) $$
and to create a pairwise dependence it is necessary to include the other weights for each other density $f_l$. This at least to us only seems possible by taking
$$w_{jk}=\sum_{l=1}^m p_{jl}w_{jlk}$$
for some weights $(p_{jl})$. And the common component is worked out by taking $w_{jlk}=w_{ljk}$.
There is no unidentifiability here because we have $w_{jlk}=w_{ljk}$ and it is this feature which is creating the dependence between the densities $(f_1,\ldots,f_m)$. 

To avoid cluttering up the notation, at this point we adopt a simpler notation
for the random densities $f_j$; namely, from now on, we
denote $f_j(x|\,({\bf w_{j1},\ldots, w_{jm})},\vartheta)$ by $f_j(x)$.

Given mutually independent observations $x=(x_{ji})$ for $j=1,\ldots,m$ and $i=1,\ldots,n_j$,
our method of inference will be the Gibbs sampler and we will rely crucially on
slice latent variables (Walker (2007), Kalli et al. (2010)); so for each $f_j$ we introduce the latent
variables $u_j=(u_{ji})_{i=1}^{n_j}$ such that the joint density of $u_{ji}$ with
$x_{ji}$ is given by
\begin{equation}
\label{fjxjuj}
f_j(x_{ji},u_{ji})\,=\,
\sum_{l=1}^m p_{jl}\sum_{k=1}^\infty {\bf 1}(u_{ji}<w_{jlk})\,K(x_{ji}|\,\theta_{k}).
\end{equation}

This augmented scheme is at the core of our sampling methodology. It essentially shifts the problem from  one of sampling from a mixture
with an infinite number of components to actually having to deal with only a finite number of them. It can be readily verified that the  sets
\begin{equation}
\label{constraintsets}
A_{w_{jl}}(u_{ji})=\{k\in{\mathbb N}:u_{ji}<w_{jlk}\},
\end{equation}
with $w_{jl}=w_{lj}$, are finite.

This means, and this will become clearer later on, that for each pair $(j,l)$ only a finite
number of the $(\theta_k)$ and
$w_{jl}=(w_{jlk})_{k=1}^\infty$, are  needed in each iteration of the Gibbs sampler.

We can then express the  $f_j$ $u_j$-augmented random densities in (\ref{fjxjuj}) as follows:
\begin{equation}
\label{fju}
f_j(x_{ji},u_{ji})\,=\,
\sum_{l=1}^m p_{jl}\sum_{k\in A_{w_{jl}}(u_{ji})}K(x_{ji}|\,\theta_{k}),\quad 1\le i \le n_j.
\end{equation}

We now introduce latent variables $\delta=(\delta_{ji})_{i=1}^{n_j}$ selecting the mixture
and $d=(d_{ji})_{i=1}^{n_j}$ selecting the component within the mixture from which the observations come; so
for each $i=1,\ldots, n_j$ we have
$$
f_j(x_{ji},u_{ji},d_{ji}|\,\delta_{ji})\,=\,\prod_{r=1}^m
                 {\bf 1}(u_{ji}<w_{jrd_{ji}})^{\delta_{ji}^r} K(x_{ji}|\,\theta_{d_{ji}})
$$
with $\delta_{ji}=(\delta_{ji}^1,\ldots,\delta_{ji}^m)$ and $\mbox{Pr}(\delta_{ji}=\hat{e}_l)=p_{jl}$,
where $\hat{e}_l$ denotes the usual basis vector having its only nonzero component equal to $1$
at position $l$.

Hence, for a sample of size $n_1$ from $f_1$, a sample of size $n_2$ from $f_2$, etc.,
a sample of size $n_m$ from $f_m$
we can write the full likelihood as a multiple product:
$$
f(x,u,d\,|\,\delta)\,=\,
\prod_{j=1}^m\prod_{i=1}^{n_j}\prod_{l=1}^m
{\bf 1}(u_{ji}<w_{jld_{ji}})^{\delta_{ji}^l}K(x_{ji}|\,\theta_{d_{ji}}).
$$
It will be reinforcing our intuition, and  will make the Gibbs sampling algorithmic steps, as well as the dependencies between variables, much
clearer if we express the model concisely in a hierarchical fashion using  the auxiliary variables and the stick breaking representation.
We thus have, for $j=1\ldots,m$ and $i=1,\ldots, n_j$,

\begin{eqnarray}
\nonumber
& & x_{ji},u_{ji}\,|\,  d_{ji},\delta_{ji},\theta_{d_{ji}},(w_{jrd_{ji}})_{r=1}^m\sim
\prod_{r=1}^m K(x_{ji}|\theta_{d_{ji}})\,\left\{{\cal U}(u_{ji}\,|0,w_{jrd_{ji}})\right\}^{\delta_{ji}^r} \\
\nonumber
& & \mbox{Pr}(d_{ji}=k\,|\,w_{ji},\delta_{ji}=\hat{e}_l)\,=\, w_{jlk}\\
\nonumber
& & \mbox{Pr}(\delta_{ji}=\hat{e}_l)\,=\,p_{jl}\\
\nonumber
& & w_{jlk}  =  z_{jlk}\prod_{r<k}(1-z_{jlr}),\,\, z_{jlk}\stackrel{\rm iid}\sim\mathrm{beta}(1,c),\,\,
\theta_{k}\,\stackrel{\rm iid}\sim\, p_0,\; k\in {\mathbb N},
\end{eqnarray}
where ${\cal U}(u|\alpha, \beta)$ is the uniform density over the interval $(\alpha, \beta)$.

It is also clear, by construction, that we can by choice of the $(p_{jl})$ and the $(g_{jl})$ arrange for $f_j$ and $f_l$ to be as close or as far apart as desired with respect to the $L_2$ metric.

\vspace{0.2in} \noindent {\bf 4. The Gibbs sampler.}
We are now ready to describe the Gibbs sampler and the full conditional densities for
estimating the model, having completed the model by assuming that the prior for
$(p_{j1},\ldots,p_{j,m-1})$ for $j=1,\ldots,m$ is Dirichlet with fixed parameters $(a_{j1},\ldots,a_{jm})$ i.e.
$$
\pi(p_{j1},\ldots,p_{j,m-1})\,\propto\,p_{j1}^{a_{j1}-1}\cdots p_{j,m-1}^{a_{j,m-1}-1}
\left(1-p_{j1}-\cdots-p_{j,m-1}\right)^{a_{jm}-1}.
$$
At each iteration we will sample variables,
\begin{eqnarray}
\nonumber
& & w_{jlk},\theta_{k},~~\, 1\le j\le l\le m,\,  1\le k\le N^*,\\
\nonumber
& & u_{ji},d_{ji},\delta_{ji},\, 1\le j\le m,\, 1\le i\le n_j,\\
\nonumber
& & p_{jl},~~~~~~~~~~ 1\le j\le m,\, 1\le l\le m-1,
\end{eqnarray}
with $N^*$ almost surely finite. In the sequel we will see how to obtain $N^*$.

{\bf A.} We start with initial $\{(d_{ji},\delta_{ji})\}$ for $j=1,\ldots,m$ and $i=1,\ldots,n_j$, and
$(p_{jl})$ for $j=1,\ldots,m$ and $l=1,\ldots,m-1$. The first task is to generate the $(w_{jlk},\theta_{k})$.
We will
do this by sampling from the conditional distribution with the $(u_{ji})$ for $j=1,\ldots,m$
and $i=1,\ldots,n_j$ integrated out. Then standard results, see Kalli et al. (2010), give
$$
\pi(z_{jjk}\,|\cdots)\,=\,{\rm beta}\left(1+\sum_{i=1}^{n_j}{\bf 1}(d_{ji}=k,\delta_{ji}=\hat{e}_j),
                                     c+\sum_{i=1}^{n_j}{\bf 1}(d_{ji}>k,\delta_{ji}=\hat{e}_j)\right),
$$
also for $j\ne l$ we have
\begin{eqnarray}
\nonumber
\pi(z_{jlk}\,|\cdots) & = & {\rm beta}\left(1+\sum_{i=1}^{n_j}{\bf 1}(d_{ji}=k,\delta_{ji}=\hat{e}_l)
                                        +\sum_{i=1}^{n_l}{\bf 1}(d_{li}=k,\delta_{li}=\hat{e}_j)\right., \\
\nonumber
                   &   & \quad\quad\; \left.c+\sum_{i=1}^{n_j}{\bf 1}(d_{ji}>k,\delta_{ji}=\hat{e}_l)
                                       +\sum_{i=1}^{n_l}{\bf 1}(d_{li}>k,\delta_{li}=\hat{e}_j)\right).
\end{eqnarray}

The $\delta$'s and $d$'s will only enter the equation for $j\leq M=\max_{i,j}\{d_{ji}\}$.
For $j>M$ we take all the $(z_{jlk})$
independently from beta$(1,c)$ and take the $(\theta_{k})$ independently from $p_0$.
The $z$'s yield the $(w_{jlk})$ according to the stick-breaking formula.

{\bf B.} Here we describe how to sample the $(\theta_{k})$ for $k\leq M$. We have
for $1\le j \le m$
\begin{equation}
\label{thetaiik}
\pi(\theta_{k}|\cdots)\,\propto\, p_0(\theta_{k})\prod_{j=1}^m
\prod_{i=1,\atop d_{ji}=k}^{n_j}K(x_{ji}|\theta_{k}).
\end{equation}

{\bf C.} Before we concern ourselves with how many of the $z$'s and $\theta$'s to sample beyond $M$,
we sample the $(u_{ji})$. From the likelihood one has
$$
\pi(u_{ji}|\cdots) \, \propto \, \prod_{l=1}^m {\bf 1}(u_{ji}<w_{jld_{ji}})^{\delta_{ji}^l},
$$
where $w_{jld_{ji}}=w_{ljd_{ji}}$ when $j>l$.
So, if $\delta_{ji}=\hat{e}_l$ we take $u_{ji}$ uniform from $(0,w_{jld_{ji}})$.
For example, when $m=3$ we have for $1\le i\le n_1$,
\begin{eqnarray}
\pi(u_{1i}|\cdots)\,=\,\left\{
\begin{array}{lll}
{\cal U}(0, w_{11d_{1i}})      & \delta_{1i}=\hat{e}_1\nonumber\\
{\cal U}(0, w_{12d_{1i}})      & \delta_{1i}=\hat{e}_2\nonumber\\
{\cal U}(0, w_{13d_{1i}})      & \delta_{1i}=\hat{e}_3,\\
\end{array}\right.
\end{eqnarray}
and so on.

{\bf D.} We now proceed to sample the rest of the $(w_{jlk})$ and $(\theta_{k})$.
Let $N_{jl}$ be the
smallest integer $N$ for which
$$
\sum_{k=1}^{N}w_{jlk}>1-u_{jl}^*,
$$
where for $j=l$ we have
$$
u_{jj}^*\,=\,\min_i\{u_{ji}\},\quad j=1,\ldots,m,
$$
also for $1\le j<l\le m$ we have
$$
u_{jl}^*\,=\,\min\{\min_i\{u_{ji}\},\min_i\{u_{li}\}\}.
$$
This then implies that we must sample, in order to sample the $(d_{ji},\delta_{ji})$, the 
rest of $(w_{jlk},\theta_{k})$ from the prior for $k=M+1,\ldots,N^*$ where $N^*=\max_{jl}\{ N_{jl}\}$.

{\bf E.}  We now concentrate on sampling the $(d_{ji},\delta_{ji})_{i=1}^{n_j}$.
To do this we first need to explicitly find the constraint sets in relations (\ref{constraintsets})
and (\ref{fju}) then for $1\le j\le m$, $1\le l\le m$ and $1\le i\le n_j$ 
the likelihood expression gives
$$
\mbox{Pr}(d_{ji}=k, \delta_{ji}=\hat{e}_l|\cdots)\,\propto\,
p_{jl}\;{\bf 1}\left( k\in A_{w_{jl}}(u_{ji})\right)K(x_{ji}|\theta_{k}),
$$
where $p_{jm}=1-\sum_{r=1}^{m-1}p_{jr}$ and $w_{jl}=w_{lj}$. Also we have
$$
\mbox{Pr}(d_{ji}=k, \delta_{ji}=\hat{e}_l|\cdots)\,=\,
{p_{jl}\;K(x_{ji}|\theta_{k})\over \sum_{s\in A_{w_{jl}}(u_{ji})}K(x_{ji}|\theta_{s})},
\quad k\in A_{w_{jl}}(u_{ji}).
$$
Now it is clear that for each $(j,i)$ we can sample the $\{d_{ji},\delta_{ji}\}$ together as a block.

{\bf F.}  It is also easily seen that the full conditional for each $(p_{j1},\ldots,p_{j,m-1})$ is a Dirichlet distribution, namely for $j=1,\ldots, m$ we have
\begin{eqnarray}
\label{DPP}
\pi(p_{j1},\ldots,p_{j,m-1}|\cdots) & \propto & p_{j1}^{\alpha_{j1}+\sum_{i=1}^{n_1}{\bf 1}(\delta_{ji}=\hat{e}_1)-1}
\cdots p_{j,m-1}^{\alpha_{j,m-1}-1+\sum_{i=1}^{n_{m-1}}{\bf 1}(\delta_{ji}=\hat{e}_{m-1})-1}\\
\nonumber
& \times & \left(1-p_{j1}-\cdots-p_{j,m-1}\right)^{\alpha_{jm}+\sum_{i=1}^{n_m}{\bf 1}(\delta_{ji}=\hat{e}_m)-1}.
\end{eqnarray}


We can use the output at each iteration of the Gibbs sampler, after a sensible burn-in
time period, to sample from the densities $(f_1,\ldots,f_m)$. So we sample  independently $(x_{j, n_{j}+1})$ from the densities based on the current parameter values of each density. This would then involve, for each $j$, sampling the component $l$ according to the probabilities $(p_{j1},\ldots,p_{j\,m-1})$ and then sampling $x_{j,n_{j}+1}$ from the $K(\cdot|\theta_{k})$ where $k$ is chosen according to the probabilities $(w_{jlk})$. These collection of samples collected over the course of the Gibbs sampler can be used to provide estimates for the $m$ densities.

\vspace{0.2in} \noindent {\bf 5. Comparing PDDP and CAPDDP priors -- Numerical Illustrations}

In this section we compare the pairwise dependent Dirichlet process (PDDP) and the common atoms 
pairwise dependent Dirichlet process (CAPDDP) models. We present  two simulated and one real data example with $m=3$.
For the choice  of kernel we have    a normal model $K(x|\theta)={\rm N}(x|\theta)$, 
where $\theta=(\mu, \lambda)$ and  $\lambda$  is the precision.
The prior for the means and the precisions for both PDDP and CAPDDP models
will be independent normals ${\rm N}\left(0,s^{-1}\right)$ and gammas $Ga(\epsilon, \epsilon)$, respectively,
i.e. $P_0(d\mu, d\lambda)\,=\,{\rm N}\left(\mu|0,s^{-1}\right)\,Ga(\lambda|\epsilon, \epsilon)\,d\mu\, d\lambda$.
Attempting a noninformative prior specification, we took in all our numerical experiments
$s=0.001$ and $\epsilon=0.001$. 
The concentration parameter is everywhere constant and has been set to $c=1$.
Our only change, across examples, will be the value 
of the hyperparameter $(\alpha_{ji})$ of the selection probabilities.

Our finding is that the massive extra cost of more parameters and the essentially equivalent predictive performance of the models combined with the lack of availability of the computation of distances, puts the PDDP model at a significant disadvantage to the CAPDDP model.

\vspace{0.1in} \noindent {\bf First simulated data example:~}
We simulated three data sets, of sizes $n_1=80$, $n_2=30$ and $n_3=80$, independently from
\begin{eqnarray}
\nonumber
x_{1i} & \sim & f_1(x)\,\,=\,\,{\rm Ga}(2-x\,|2,1),\quad 1\le i\le n_1\\
\label{simdataeg1}
x_{2i} & \sim & f_2(x)\,\,=\,\,{\rm N}(x\,|0,2),\quad 1\le i\le n_2\\
\nonumber
x_{3i} & \sim & f_3(x)\,\,=\,\,{\rm Ga}(x+2\,|2,1),\quad 1\le i\le n_3.
\end{eqnarray}

In both cases, when atoms are common (CAPDDP) and when atoms are unequal (PDDP),  
the hyperparameters in relation (\ref{DPP}) of the Dirichlet priors for the mixing weights are given by
$\alpha_{ji}={\bf 1}(j\ne i)+3\:{\bf 1}(j=i)$. We sample $70,000$ points from the predictives
after a burn-in of $10,000$.

The following results are presented in Figures 1,2,3,4:
\begin{enumerate}
\item  In  $1$(a), $1$(b) and  $1$(c) the true densities, as well as the  kernel density estimates based on samples of the first $40,000$ points after the  burn-in period for  the PDDP and CAPDDP models, are superimposed over the corresponding data 
sets. We note the similarity of the posterior predictive estimates of 
the densities $f_1$, $f_2$ and $f_3$. 
 \item In  $1$(d), $1$(e) and  $1$(f)  the histograms of the predictive samples 
coming from the PDDP model along with the associated KDE curves.
\item In  $1$(g), $1$(h) and  $1$(i) the corresponding predictive samples and KDE curves of the CAPDDP model. 

\item In $2$(a), $2$(b) and $2$(c)  the histograms of the sampled values 
of the conditonal expectations of the $L_2$ distances via
\begin{equation}
\label{dapprox}
\E\left[\, d_2(f_j,f_i)\, |{\bf w_j}, {\bf w_i}, N^* \right]\,\propto\,
\sum_{k=1}^{N^*} (w_{jk}-w_{ik})^2,\quad 1\le j<i\le 3,
\end{equation}
where
\begin{eqnarray}
\nonumber
w_{1k} & = & p_{11}w_{11k}+p_{12}w_{12k}+(1-p_{11}-p_{12})\,w_{13k}\\
\label{wdef}
w_{2k} & = & p_{21}w_{12k}+p_{22}w_{22k}+(1-p_{21}-p_{22})\,w_{23k}\\
\nonumber
w_{3k} & = & p_{31}w_{13k}+p_{32}w_{23k}+(1-p_{31}-p_{32})\,w_{33k}.
\end{eqnarray}

\item In  $2$(d), $2$(e) and $2$(f) the running averages of
the associated expected $L_2$ distances 
$$
d_{ji}=\E\left[\, d_2(f_j,f_i)\right],\quad 1\le j<i\le 3.
$$
Our numerical approximations after $50,000$ iterations are $d_{12}\propto 0.0969$, $d_{23}\propto 0.0970$ and  
$d_{13}\propto 0.2555$. These values exhibit the same features as the $L_2$ distances
of the true densities given in relations (\ref{simdataeg1})
\begin{eqnarray}
\nonumber
d_2(f_1,f_2) & = & d_2(f_2,f_3)\,\, =\,\,
{1\over 4}+{1\over 2\sqrt{2\pi}}-{2\over e\sqrt{\pi}}\,\,\approx \,\, 0.0346\\
\nonumber
d_2(f_1,f_3) & = & {64\over 3}\,e^{-4}\,\,\approx\,\, 0.1093.
\end{eqnarray} 
\item In $3$(a), $3$(b) and $3$(c) the histograms of the $p$-values of the Anderson--Darling two-sample test between 
$700$ samples of size $100$ coming from the PDDP and CAPDDP models. The test rejects $267$ samples out of 
$700$ of the $f_1$-predictive, rejects $234$ samples out of $700$ of the $f_2$-predictive 
and $263$ samples out of $700$ of the $f_3$-predictive.

\item In $4$(a), $4$(b) the histograms of the $p$-values of one-sample Anderson--Darling 
tests. In $4$(a) we test $700$ samples of size $100$ against the hypothesis that
the samples from the $f_2$-PDDP-predictive are coming from a normal with mean $0$ and variance $2$. 
The rejection rate is $242/700$. The rejection rate is smaller  for the corresponding $f_2$-CAPDDP-predictive,
namely $216/700$.
\end{enumerate}

\vspace{0.1in} \noindent {\bf Second simulated data example:~} We have simulated idependently
two groups of data sets, the ``large" group 
$$
G_1=\left\{(x_{1i_1})_{i_1=1}^{300},(x_{2i_2})_{i_2=1}^{300},(x_{3i_3})_{i_3=1}^{300}\right\},
$$ 
and the ``small" group
$$
G_2=\left\{(x_{1i_1})_{i_1=1}^{120},(x_{2i_2})_{i_2=1}^{60},(x_{3i_3})_{i_3=1}^{120}\right\},
$$ 
from the normal $3$-mixtures
\begin{eqnarray}
\nonumber
x_{1i_1} & \sim & f_1(x)\,\,=\,\,{1\over 3}\,{\rm N}(x|-10, 1)+{1\over 3}\,{\rm N}(x|-20, 1)+{1\over 3}\,{\rm N}(x|20, 1)\\
\label{simdataeg2}
x_{2i_2} & \sim & f_2(x)\,\,=\,\,{1\over 3}\,{\rm N}(x|-20, 1)+{1\over 3}\,{\rm N}(x|  0, 1)+{1\over 3}\,{\rm N}(x|30, 1)\\
\nonumber
x_{2i_3} & \sim & f_3(x)\,\,=\,\,{1\over 3}\,{\rm N}(x| 20, 1)+{1\over 3}\,{\rm N}(x| 30, 1)+{1\over 3}\,{\rm N}(x|10, 1).
\end{eqnarray}
In both cases, CAPDDP and PDDP,  the hyperparameters in relation (\ref{DPP}) 
of the Dirichlet priors for the mixing weights $p_{ji}$ are 
$\alpha_{ji}=1$. We sample $70,000$ points from the predictives
after a burn-in of $10,000$.

The following results are presented in Figure 5:

\begin{enumerate}
\item In  $5(a)$--$5(c)$ the $f_j$-CAPDDP-predictives for the large group of data sets $G_1$. As it
can be seen they predict very effectively the true densities apart from the nearly discernible
small bumps in the areas around the modes of the $3$-mixtures. Such mass dislocations can occur as well during the 
computation of the $f_j$-PDDP-predictives but they quicly dissappear as the MCMC evolves. In the common atoms 
case such  mass dislocations need a very large number of iterations to completely dissappear.
\item In $5(d)$--$5(f)$  the KDE curves of the 
distributions of the conditonal expectations of the $L_2$ distances. As it can be seen, due to the large sample size,
these are tight and centered about the same value $d_{ji}\approx 0.44$.
\item In  $5(g)$--$5(i)$ the $f_j$-CAPDDP-predictives for the small group of data sets $G_2$, which is a reduction of the $G_1$ group. Here the mass disslocations expand due to the small sample sizes.
Nevertheless all the $L_2$ conditonal expectation distributions, remain centered to $d_{ji}\approx 0.4$.
We observe that the variance has been increased considerably purely due to the smaller sample sizes.
Note that the $L_2$ distances of the true densities given in relations (\ref{simdataeg2}) are  
$$
d_2(f_1,f_2)\, =\,d_2(f_1,f_3)\, =\,d_2(f_2,f_3)\,\approx\, 0.125.
$$
\end{enumerate}

\vspace{0.1in} \noindent {\bf Real data example:~}
The data set is to be found at 
$http://lib.stat.cmu.edu/datasets/pbcseq$ 
and involves data from 312 individuals. We take the observation as SGOT (serum glutamic-oxaloacetic transaminase) level, just prior to liver transplant or death or the last observation recorded, under three conditions on the individual
\begin{enumerate}
\item{The individual is dead without transplantation.}
\item{The individual had a transplant.}
\item{The individual is alive without transplantation.}
\end{enumerate}
We normalize the means of all three data sets to zero. Since it is reasonable to assume the densities for the observations are similar for the three categories (especially for the last two), we adopt the model proposed in this paper with $m=3$. The number of transplanted individuals is small (sample size of 28) so it is reasonable to borrow strength for this density from the other two.  
We took the hyperparameters of the Dirichlet priors for the mixing weights $(p_{ji})$ are for $1\le i,j\le 3$
\[ 
\alpha_{ji} = \left\{ 
\begin{array}{l l}
10 & \quad \mbox{if $j=i=1$ or $j=i=3$}\\
1 & \quad \mbox{elsewhere}.\\
\end{array} \right. 
\]
We sample $50,000$ points from the predictives after a burn-in of $10,000$.

The following results are presented in Figure 6:
\begin{enumerate}
\item In $6$(a)--$6$(c).The histograms for the three data sets  
along with the superimposed    predictive density curves of the $f_j\,$-PDDP (solid) and  
$f_j\,$-CAPDDP(dashed). 
\item In  $6$(d)--$6$(f)  the approximate
densities of the $50,000$ sampled values of the conditonal expectations of the $L_2$ distances.

We note how the distribution of $\E\left[\, d_2(f_2,f_3)\, |{\bf w_2, w_3} , N^* \right]$ concentrates
near zero, and the similarity of the distributions of 
$\E\left[\, d_2(f_1,f_2)\, |{\bf w_1, w_2}, N^* \right]$
and $\E\left[\, d_2(f_1,f_3)\, |{\bf w_1, w_3}, N^* \right]$. 

\item In  $6$(g)--$6$(i) the running averages.
\end{enumerate}

\vspace{0.1in} \noindent 
General distinctive features between the CAPDDP and the PDDP models include: 
\begin{enumerate}
\item{The average number of clusters during CAPDDP computations is larger than the average number 
of clusters in PDDP computations. A detailed comparison between the numerical approximations of the average
number of clusters, for both the CAPDDP and the PDDP case, can be found in Table $1$ where we compare the 
average number of clusters coming from the three numerical examples.}
\item{Borrowing of strength, seems to be sometimes weaker for the undersampled data set in the CAPDDP case. 
A detailed comparison on the borrowing of strength for the second predictiive in the three
numerical examples can be found in Table $2$. The larger the predictive running average for the $p_{22}$ 
selection probability is, the weaker the borrowing of strength between the second data set and the other
two covariate data sets.} 
\end{enumerate}  

\vspace{0.2in} \noindent {\bf 6. Discussion.} We have shown when constructing pairwise dependent random densities, 
for the purposes of borrowing strength, the random Dirichlet processes used to effect this can be taken to have 
identical atoms. Random probability measures can be constructed by adapting the weights to these atoms and can provide 
arbitrary degrees of proximity; from identical to far apart. These distances can be readily evaluated precisely because 
they have identical atoms and hence share the same supports.
Also the time complexity difference between the two algorithms is substancial. Let $\Delta T$ be the time complexity diffrerence
between the PDDP and the CAPDDP algorithms based on a single sweep of the associated Gibbs samplers. Then it is not difficult to 
verify that
$$
\Delta T \propto \left({m(m+1)\over 2}-1\right) N^*\sum_{j=1}^m n_j,
$$
where $N^*$ is a Poisson random variable with mean $c\log(1/u^*)$ and $u^*$ being the global minimum of the auxiliary variables $u_{ji}$'s for $1\le j\le m$
and $1\le i\le n_j$ (Muliere, Tardella 1998). This time complexity difference is due to the fact that in the uncommon atom case we have to sample the posterior 
locations of $m(m+1)/2$ different random measures.

We can also extend this principle to regression scenarios whereby densities $(f_j)$ are characterized by covariates $z$. We can still adopt the idea of common atoms by designing
$$f_z(x)=\int K(x|\theta)\mP_z(d\theta)$$
where
$$\mP_z(d\theta)=\sum_{k=1}^\infty w_k(z)\delta_{\theta_k}(d\theta).$$
So as before we would have
$$f_z(x)=\sum_{k=1}^\infty w_k(z)K(x|\theta_k)$$
and the covariate dependent weights will adapt the weights on the atoms to change the shape of the density. This could be densities which are far apart, based on a covariate pair which are far apart, or densities close to each other when the covariates are close. Indeed this is a simpler version of many Bayesian nonparametric regression models which typcially include $z$ in $K(x|z,\theta)$ and have $\theta_k(z)$. We believe these extras are unnecessary.    

To conform with our model described in Section 3 we could now take
$$w_k(z)=\int p(s,z)w_{k}(s,z)ds$$
with $w_k(z,s)=w_k(s,z)$ and $\int p(s,z)ds=1$ for all $z$.
So, for each $(s,z)$, for $s\geq z$, $(w_k(s,z))$ is a stick-breaking set of weights, mutually independent over $(s\geq z)$, and
for each $z$, $(p(\cdot,z))$ is a Dirichlet process, mutually independent over $(z)$.    

This is not as complicated as it looks, essentially we would replace a finite $m$ with $m=+\infty$. So if covariates are from the set $(z_1,z_2,z_3,\ldots)$ then we would model
$$f_j(x)=f_{z_j}(x)=\sum_{k=1}^\infty \sum_{l=1}^\infty p_{jl}w_{jlk}K(x|\theta_k).$$
Of course, in practice, only a finite number of covariates would be observed.

\vspace{0.2in} \noindent {\bf Acknowledgement.} The work was completed during a visit by the third author to the University of the Aegean in Karlovassi, Samos.

\newpage

\begin{table}[h]
\begin{center}
\caption{Running average for the number of clusters.}
\label{tab:table1}
\begin{tabular}{|l|l|c|c|c|}
\hline
Data sets and sample sizes       & Model   & Predictive 1 &  Predictive 2 &  Predictive 3  \\ \hline\hline
 gamma-normal-gamma            & CAPDDP  &   3.075      &   4.456       &   4.044        \\
 $n_1=100~n_2=50~n_3=100$        & PDDP    &   1.571      &   1.079       &   2.183        \\ \hline
 normal 3-mixtures              & CAPDDP  &   4.861      &   4.773       &   3.154        \\
 $n_1=120~n_2=60~n_3=120$        & PDDP    &   3.090      &   3.041       &   3.062        \\ \hline
 real example                    & CAPDDP  &   5.045      &   6.178       &   4.328        \\
 $n_1=143~n_2=28~n_3=139$        & PDDP    &   3.090      &   3.041       &   3.062        \\ \hline           
\end{tabular}
\end{center}
\end{table}

\begin{table}[h]
\begin{center}
\caption{Running averages for the selection probabilities in predictive $f_2$.}
\label{tab:table2}
\begin{tabular}{|l|l|c|c|c|}
\hline
Data sets and sample sizes       & Model   &   $p_{21}$   &  $p_{22}$     &  $p_{23}$     \\ \hline\hline
 gamma-normal-gamma            & CAPDDP  &   0.343      &   0.301       &   0.356        \\
 $n_1=100~n_2=50~n_3=100$        & PDDP    &   0.484      &   0.167       &   0.349       \\ \hline
 normal 3-mixtures              & CAPDDP  &   0.146      &   0.702       &   0.152        \\
 $n_1=120~n_2=60~n_3=120$        & PDDP    &   0.359      &   0.405       &   0.236        \\ \hline
 real example                    & CAPDDP  &   0.341      &   0.325       &   0.334       \\
 $n_1=143~n_2=28~n_3=139$        & PDDP    &   0.305      &   0.328       &   0.366       \\ \hline           
\end{tabular}
\end{center}
\end{table}

\newpage

\vspace{0.2in} \noindent {\bf References.}

\begin{description}

\item Bulla, P., Muliere, P. and Walker, S.G. (2009). A Bayesian nonparametric
estimator of a multivariate survival function. {\sl Journal of Statistical
Planning and Inference} {\bf 139}, 3639--3648.

\item De Iorio, M., M\"uller, P., Rosner, G.L. and MacEachern, S.N. (2004). An ANOVA model for dependent random measures. {\sl Journal of the American Statistical Association}
{\bf 99}, 205--215. 

\item Dunson, D.B. and Park, J.H. (2008). Kernel stick-breaking processes. {\sl Biometrika}   {\bf 95},   307--323.

\item Ferguson, T.S. (1973). A Bayesian analysis of some nonparametric problems.
{\sl Annals of Statistics} {\bf 1}, 209-230.

\item Griffin, J.E. and Steel, M.F.J. (2006). Order-based dependent Dirichlet processes. {\sl Journal of the American Statistical Association}
{\bf 101}, 179--194.

\item Griffin, J.E., Kolossiatis, M. and Steel, M.F.J. (2013). Comparing distributions by using dependent normalized random-measure mixtures. {\sl Journal of the Royal Statistical Society, Series B} {\bf 75}, 499--529.

\item Hatjispyros, S.J., Nicoleris, T. and Walker, S.G. (2011). 
Dependent mixtures of Dirichlet processes.  {\sl Computational Statistics and Data Analysis} {\bf 55}, 2011--2025.

\item Hjort, N.L., Holmes, C.C., M\"uller, P. and Walker, S.G. (2010). {\sl Bayesian Nonparametrics}. Cambridge University Press. 

\item Ishwaran, H. and James, L.F. (2001). Gibbs sampling methods for stick-breaking priors.
{\sl Journal of the American Statistical Association} {\bf 96}, 161--173.

\item Kalli, M., Griffin, J.E. and Walker, S.G. (2010). Slice Sampling Mixture Models. {\sl Statistics and Computing}
{\bf 21}, 93--105.

\item Kolossiatis, M., Griffin, J.E. and Steel, M.F.J. (2013). On Bayesian nonparametric modelling of two correlated distributions. {\sl Statistics and Computing} {\bf 23}, 1--15.  

\item Lo, A.Y. (1984). On a class of Bayesian nonparametric estimates I. Density estimates.
{\sl Annals of Statistics} {\bf 12}, 351--357.

\item MacEachern, S.N. (1999). Dependent nonparametric processes. In {\sl ``Proceedings of the Section on Bayesian Statistical Science"}, pp.
50--55. American Statistical Association.

\item Muliere, P.,  Tardella, L. (1998). Approximating distributions of random
functionals of Ferguson-Dirichlet priors. {\sl Canadian Journal of Statistics} {\bf 26}, 283--297.

\item M\"uller, P., Quintana, F. and Rosner, G. (2004). A method for combining inference across related nonparametric Bayesian models.
{\sl Journal of the Royal Statistical Society, Series B} {\bf 66}, 735--749.

\item Sethuraman, J. (1994). A constructive definition of Dirichlet priors. {\sl Statistica Sinica} {\bf 4}, 639--650.

\item Walker, S.G. (2007). Sampling the Dirichlet mixture model with slices. {\sl Communications in Statistics} {\bf 36}, 45--54.

\end{description}

\newpage

\begin{figure}
\centerline{\includegraphics[scale=0.8]{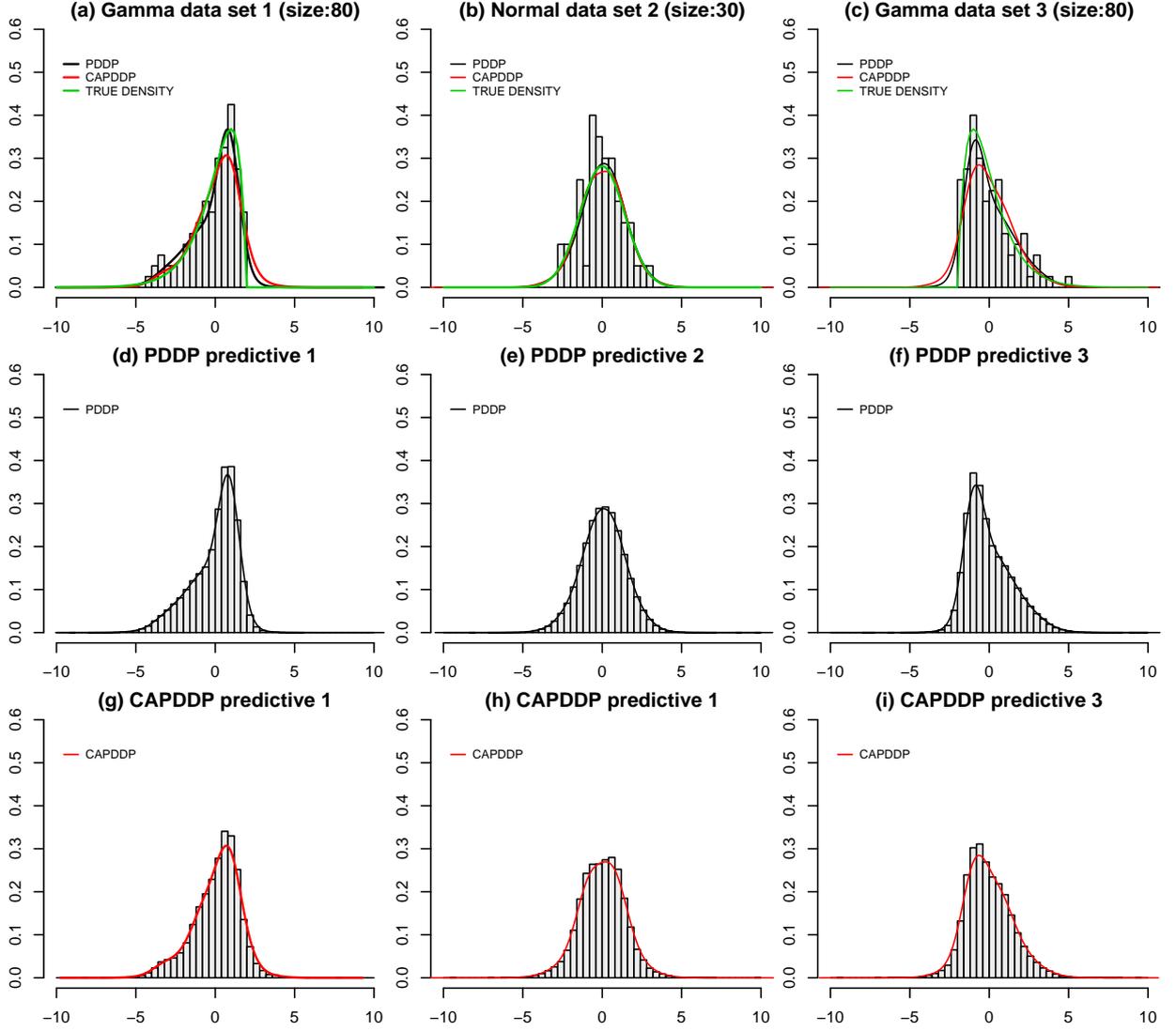}}
\caption{  (a), (b), (c): histograms for the three data sets.
           (d), (e), (f): predictives for $f_1$, $f_2$ and 
          $f_3$ in the case of the uncommon atoms model PDDP. 
          (g), (h), (i): the corresponding predictives for the case
          of the common atoms model CAPDDP. In both cases the prior 
          specifications are the same, $c=1$, $\epsilon=0.001$ and $\tau=0.001$. The 
          hyperparameters in relation (\ref{DPP}) of the Dirichlet priors are equal to 
          $\alpha_{ji}={\bf 1}(j\ne i)+3\:{\bf 1}(j=i)$.}
\end{figure}
\begin{figure}
\centerline{\includegraphics[scale=0.8]{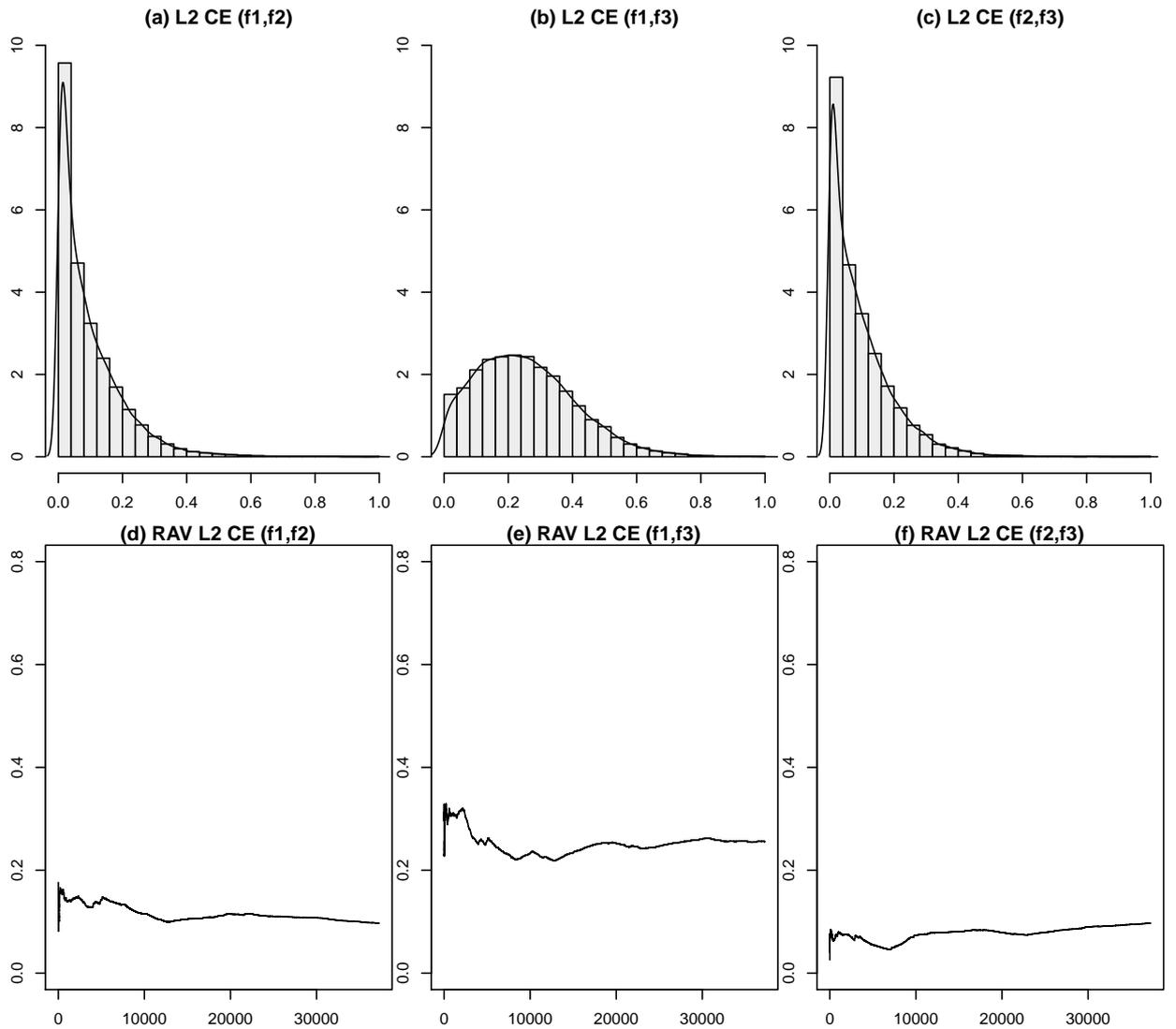}}
\caption{  Histograms
         of  the  conditional expectations a), b), c)  and   of the running averages d), e), f)
         of the $L_2$ distances from the CAPDDP model. The same prior specifications apply as in Figure $1$.}
\end{figure}

\begin{figure}
\centerline{\includegraphics[scale=0.8]{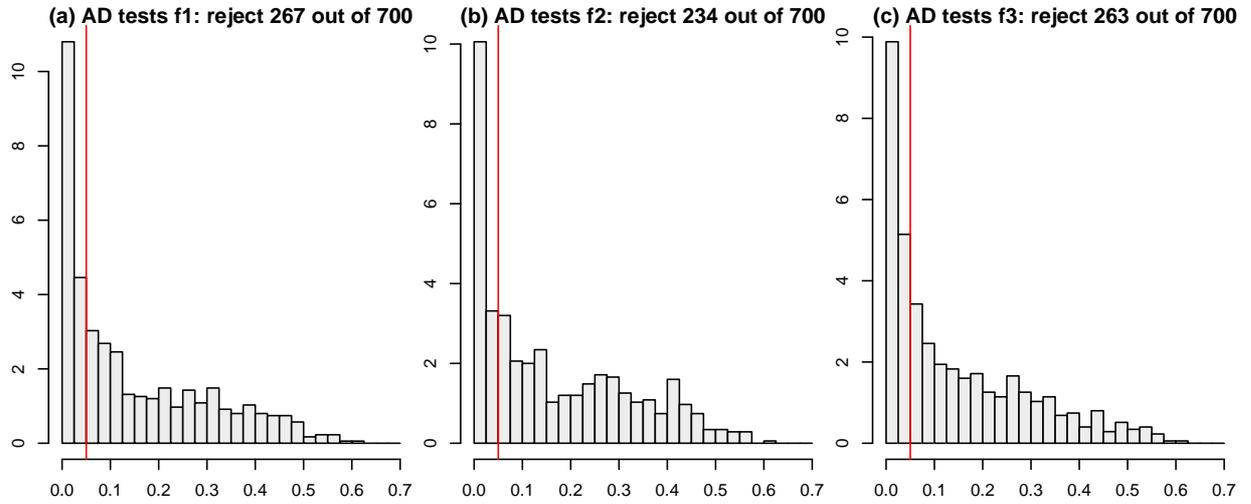}}
\caption{Histograms of p-values
         for $2$-sample Anderson--Darling tests. Test  of $700$ pairs of predictive samples 
         of size $100$, coming from the PDDP and CAPDDP models, respectively.
         Burn-in period of $30,000$.}
\end{figure}

\begin{figure}
\centerline{\includegraphics[scale=0.8]{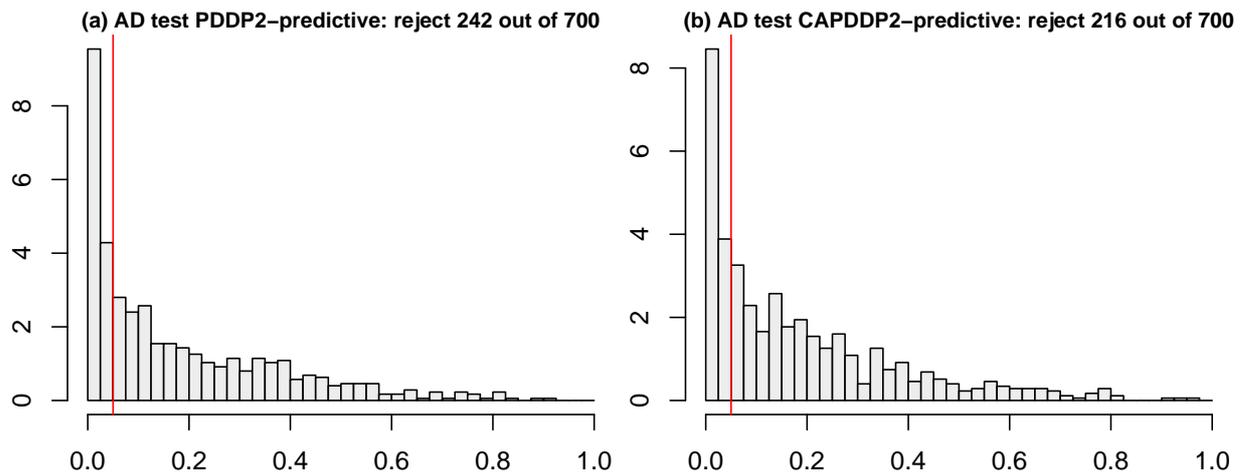}}
\caption{Histograms of p-values from $1$-sample Anderson--Darling tests, 
         between $f_2$--predictive and $N(0,2)$. 
        (a), (b): test of  $700$ pairs of predictive samples of size $100$, 
         coming from the PDDP and CAPDDP model, respectively.
         Burn-in period of $10,000$.}
\end{figure}

\begin{figure}
\centerline{\includegraphics[scale=0.8]{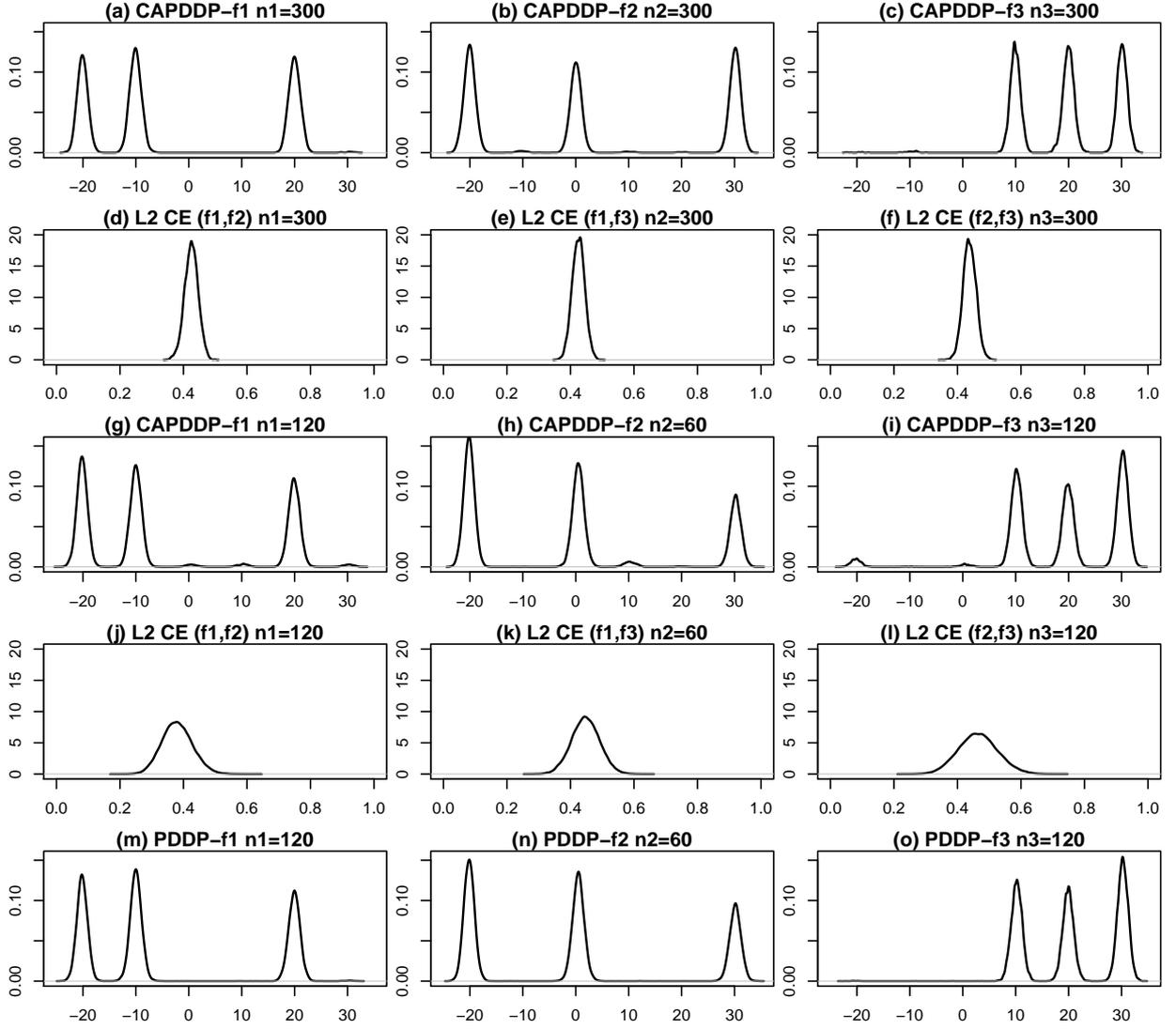}}
\caption{Predictives for $f_1$, $f_2$ and $f_3$ coming from the CAPDDP model and expected
         $L_2$ distances, for different sample sizes.
         (a)--(f): equal sample sizes $n_1=n_2=n_3=300$.
         (g)--(l): the initial sample is reduced  to size $n_1=n_3=120$ and $n_2=60$.
         In both cases the prior specifications are the same, $c=0.1$, 
         $\epsilon=0.001$ and $\tau=0.001$. The 
         hyperparameters in relation (\ref{DPP}) of the Dirichlet priors are equal to 
         $\alpha_{ji}=1$. Burn-in period of $10,000$.}
\end{figure}

\begin{figure}
\centerline{\includegraphics[scale=0.8]{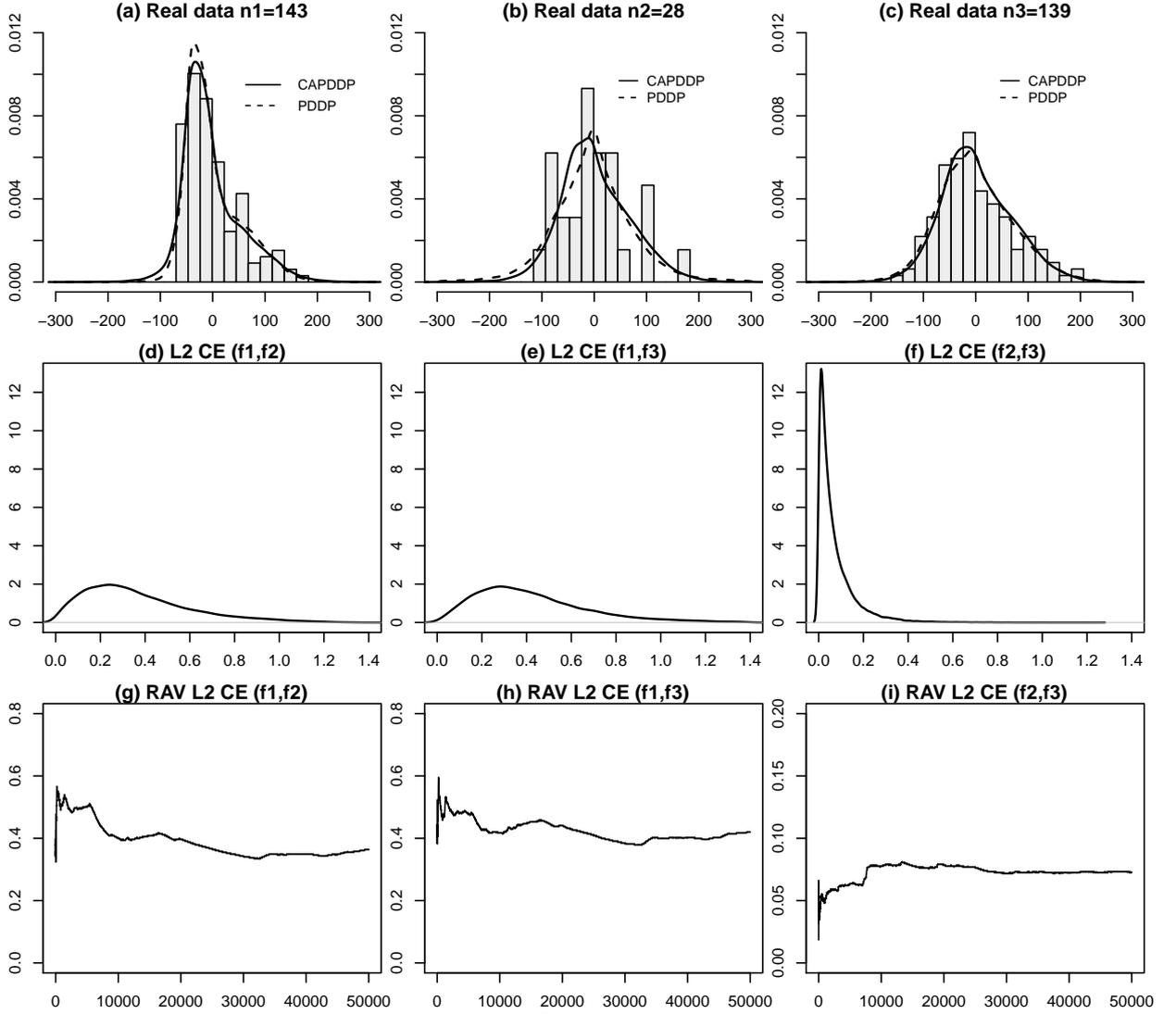}}
\caption{(a), (b), (c): histograms for the three real data sets  
         with the CAPDDP and PDDP density estimations superimposed.
          (d), (e), (f): the associated  expected 
         pairwise $L_2$ distances $L_2$ $CE (f_i,f_j)$. 
         (g), (h) (i): the running 
         averages corresponding to the expected pairwise distances RAV $L_2$ CE $(f_i,f_j)$.
         In all cases the prior 
         specifications are the same, $c=1$, $\epsilon=0.001$ and $\tau=0.001$. The 
         hyperparameters in relation (\ref{DPP}) of the Dirichlet priors are  $\alpha_{ji}=1$,
         for all $1\le j,i\le 3$. Burn-in period of $10,000$.}
\end{figure}

\end{document}